\documentclass[12pt]{amsart}
\usepackage{psfig}
\let\cal\mathcal

\title{Bases in Systems of Simplices and Chambers}

\author{Tatiana Alekseyevskaya}

\thanks{Supported by the Gabriella and
Paul Rosenbaum Foundation; part of these results was obtained at
MSRI supported by NSF grant DMS 9022140.} 

\address{Department of Mathematics\\
Rutgers University\\
Piscataway, NJ 08855\\
USA}

\email{tva@math.rutgers.edu}

\date{\today}
\newtheorem{definition}{Definition}[section]
\newtheorem{theorem}[definition]{Theorem}

\newtheorem{lemma}[definition]{Lemma}

\newtheorem{proposition}[definition]{Proposition}

\begin{document}

\begin{abstract}
We consider a finite set $E$ of points in the $n$-dimensional affine
space and two sets of objects that are generated by the set $E$: the
system $\Sigma$ of $n$-dimensional simplices with vertices in $E$ and
the system $\Gamma$ of chambers. The incidence matrix $A= \parallel
a_{\sigma, \gamma}\parallel$, $\sigma \in \Sigma,                                              \ \gamma \in \Gamma$,
induces the notion of linear independence among simplices (and among
chambers). We  present an algorithm of construction of bases of
simplices (and bases of chambers). For the case $n=2$ such an algorithm
was described in \cite{A}. However, the case of $n$-dimensional space
required a different technique. It is also proved that the constructed
bases of simplices are geometrical (according to \cite{A}).

\end{abstract}

\maketitle
\section{Introduction.}

Let $E=(e_1,e_2, \ldots , e_N), \ N>n, $  be a finite set of points in
an $n$-dimensional affine space $V$. Let $P= conv(E)$ be the convex hull
of $E$.
Let $\sigma=\sigma(e_{i_1}, \ldots , e_{i_{n+1}}) $ be the
$n$-dimensional simplex with the vertices $e_{i_1}, \ldots , e_{i_{n+1}}
\in E$. Denote by $\Sigma$ the set of all such simplices $\sigma$. All
the simplices  $\sigma$ (as a rule overlapping)  cover the polytope $P$.
The simplices $\sigma $ divide the polytope $P$ into a finite number of
chambers $\gamma$ (see Definition (\ref{def})). Denote by $\Gamma$ the
set of all chambers $\gamma$ in $P$.

\begin{definition}
\label{def}
Let $\sigma \in \Sigma$ and $\tilde{\sigma}$ be the boundary of
$\sigma$. Let $\tilde{\Sigma}= \bigcup_{\sigma \in \Sigma}
\tilde{\sigma}$ and $\bar{P}= P \setminus \tilde{\Sigma}$. Let
$\bar{\gamma}$ be a connected component of $\bar{P}$ and $\gamma$ be
closure of $\bar{\gamma}$.
We call $\gamma$ a chamber and $\bar{\gamma}$ an open chamber.
\end{definition}

Let $A$ be the incidence matrix between simplices and chambers, i.e.
$$\parallel a_{\sigma,\gamma}\parallel =1 \hbox{ iff }
\gamma \subset \sigma.$$

Consider the linear space $V_\Sigma$ generated by the rows of $A$
and the linear space $V_\Gamma$ generated by the columns of $A$
over some field of characteristic 0. Due to one-to-one correspondence
between the rows of $A$ and the simplices $\sigma \in \Sigma$,
we can speak about a linear combination of simplices instead of a linear
combination of the corresponding rows of $A$. An important question is
to construct  bases of simplices and bases of chambers, i.e. bases in
$V_\Sigma$ (or in $V_\Gamma$) that consists of simplices (or chambers)
and not of their linear combinations.

\smallskip

A basis of simplices can be also defined as follows.
Let $\phi_{\sigma} (x)$ be the characteristic function of a
simplex $\sigma$, i.e. 
$$\phi_ \sigma (x) =1,\ \if\ \ x \in \sigma\ \  and\ \  \phi_\sigma
(x) =0,\  \  x \not \in \sigma.$$
A basis of simplices is a maximal subset of simplices such that their
characteristic functions $\phi _\sigma (x)$ are linearly independent.

\smallskip

In this paper we will describe (Section 2) the inductive algorithm of
constructing bases of simplices and bases of chambers in the
$n$-dimensional affine space; the algorithm uses the case $n=2$ (see
\cite{A}) as the first step of induction. In Sections 3, 4 we prove that
the set $B$ of simplices and the set $B'$ of chambers constructed by the
algorithm are indeed bases in $V_\Sigma$ and in $V_\Gamma$ respectively.

\section{Construction of a basis of simplices and a
basis of chambers.}

\subsection{A special ordering of points $e \in E$ and related
polytopes.}

Let $E =\{e_1, \ldots , e_N\}$ be a set of points in an $n$-dimensional
affine space $V^n$, $\Sigma$ the set of $n$-dimensional simplices
$\sigma$ with the vertices in $E$ and $\Gamma$ the set of chambers
$\gamma$ (defined in Introduction.)

We will define an ordering of points $e_i \in E$ which is essential
in the construction.

\begin{lemma}
\label{lemma}
Let $E=\{e_1,e_2, \ldots, e_N \}$ be a finite set of points in the
$n$-dimensional affine space. There exists an ordering
$e_{i_1}, e_{i_2},\ldots, e_{i_N}, \ \ e_{i_k}\in E, $
such that $for \ k=1, \ldots, N$
\begin{equation}
\label{conv}
conv(e_{i_1}, \ldots , e_{i_k}) \cap conv(e_{i_{k+1}}, \ldots , e_{i_N})
=\emptyset,
\end{equation}
and there exists a hyperplane $H_k$ which separates the polytopes
\begin{equation}
\label{setF}
F_k =conv (e_{i_1}, \ldots , e_{i_k})
\end{equation}
and
\begin{equation}
\label{setE}
P_k =conv(e_{i_{k+1}},\ldots , e_{i_N})=
conv (E \setminus (e_{i_1}, \ldots , e_{i_k})).
\end{equation}
We assume also
$F_0 =\emptyset , P_N = \emptyset, P_0 =conv (E)=P, F_N =conv (E)=P$.
\end{lemma}

This lemma is proved, for example, in \cite{A}. The ordering $e_1,
\ldots, e_N$ satisfying Lemma \ref{lemma} yields a shelling of the
polytope $P$, see \cite{DK}. \smallskip

Let the ordering $e_1, \ldots, e_N$ satisfy
(\ref{conv}). Consider the sequence of polytopes $P_0=P, P_1, \ldots,
P_N=\emptyset $ defined by formula (\ref{setE}). We have $P_0 \supset
P_1 \supset \ldots \supset P_N$.
Let us denote
\begin{equation}
\label{SK}
S_k=\overline{P_{k-1} \setminus P_k},
\end{equation}
where $\overline{A}$ means the closure of the set $A$.
Let $int(S)$ be the interior of $S$.

It is easy to check that the following statements are true:
\begin{enumerate}
\item{ Each ordering $e_1, \ldots , e_N$
determines a decomposition of the polytope $P$:
\begin{equation}
\label{decomp}
P =\bigcup_{i=1}^{N-n}S_i,
\end{equation}
where $int(S_i)\cap int(S_j) = \emptyset$ for $i \neq j$.}
\item{$int(S_k) \cap E = \emptyset$ for $k=1, \ldots, N-n$. }
\item{The polytope $S_k$ is part of the convex cone with the vertex
$e_k$ and bounded by some part $\cal{L}_k$ of the boundary of $P_{k-1}$,
(where $\cal{L}_k= S_k \cap P_{k-1}$).
Note that the polytope $S_k$ is not necessarily convex\footnote{The
polytopes $S_k$ are considered in more detail in \cite{A}.}.}
\item{\begin{equation}
\label{SP}
P= S_1 \cup \ldots \cup S_{k-1} \cup P_{k-1}
\end{equation}}
\end{enumerate}

\subsection{ A map from the polytope $P$ to the hyperplane
$H_k$.}

Let an ordering of points $e_1, \ldots , e_N$ satisfy
(\ref{conv}). Consider the point $e_k \in P_{k-1}$ and the corresponding
hyperplane $H_k$. For a segment $(e_k, e_i), \ i=k+1, \ldots, N$ let
us denote
\begin{equation}
\label{EIK}
e^k_i=(e_k, e_i) \cap H_k.
\end{equation}

Thus, on the hyperplane $H_k$ we obtained the set $E_k$ of
points $e^k_1, e^k_2, \ldots , e^k _{N_k}$. (Note that $N_k \leq
N-k$ since the points $e \in E$ are not necessarily in general
position.)
In the hyperplane $H_k$ we use the points $e^k \in E_k$ to
construct simplices $\sigma^k$ and chambers $\gamma^k$ in the same way
as it was done
for the set $E$ in $V^n$ . Let $\Sigma_k$ be the set of all these
simplices $\sigma^k$ and $\Gamma _k$ the set of all chambers $\gamma^k$.

Let  $\sigma^k=\sigma(e^k_{i_1}, \ldots , e^k_{i_n}) \in \Sigma_k$.
We denote by $\overrightarrow{(e_k, e^k_{i_j})}$ the ray
starting at $e_k$ and passing through the point $(e^k_{i_j})$.
Consider the following map:
\begin{equation}
\label{mu}
\mu :\sigma^k \mapsto \sigma,\ \
where \ \sigma= \sigma(e_k, e_{i_1}, \ldots ,e_{i_n}),
\end{equation}
and where $e_{i_j}\in E$ is the nearest point to the
point $e_k$ on the ray $\overrightarrow{(e_k, e^k_{i_j})}$.
It is clear that $\sigma= \mu(\sigma^k) \in \Sigma$. The map $\mu$ is an
injection and has the following easy to check property.

\begin{proposition}
\label{prop1}
Let  $\sigma^k, \sigma^k_0 \in \Sigma_k$ be open
simplices such that $\sigma^k \cap \sigma^k_0 = \emptyset$. Then
$\mu (\sigma^k ) \cap \mu (\sigma^k_0) =\emptyset$.
\end{proposition}

Consider the set of points $E_k$ in the hyperplane $H_k$.
Let us reorder $e^k \in E_k$ according to (\ref{conv}). Let
$e^k_1, e^k_2, \ldots , e^k _{N_k}$ be such an ordering.
Similarly to formulas (\ref{setE}) and (\ref{SK}) we denote

\begin{equation}
\label{PS}
P^k=P^k_0= conv(E_k),\ P^k_j= conv(e^k_{j+1}, \ldots , e^k_{N_k}), \ \
S^k_j= P^k_{j-1} \setminus P^k_j
\end{equation}
where $j= 1, \ldots, N_k$. For these $(n-1)$-dimensional
polytopes the formulas analogous to (\ref{decomp}) and (\ref{SP}) hold:
\begin{equation}
\label{PK}
P^k =\bigcup_j S^k_j,
\end{equation}
where $int(S^k_i) \cap int( S^k_j)= \emptyset$ for $i \neq j,$
and
\begin{equation}
\label{PK1}
P^k= S^k_1 \cup \ldots \cup S^k_{j-1} \cup P^k_{j-1}.
\end{equation}

On the polytopes $S^k_j$ and $P^k_j$ let us define the following
map $\tau$: $\tau(S^k_j)$ is an $n$-dimensional cone with the vertex
$e_k$ and generated by the rays $\overrightarrow{(e_k, x)}$, where $x
\in S^k_j$. The map $\tau(P^k)$ is defined similarly. Clearly, the
following decomposition holds:

\begin{equation}
\label{tau}
\tau(P^k) = \bigcup _i \tau(S^k_i),
\end{equation}

where $int(\tau(S^k_i)) \cap int( \tau(S^k_j))= \emptyset$ for $i \neq j$.

\subsection{ Algorithm of construction of the set $B$ of simplices and
the set $B'$ of chambers.}

We will construct the set $B \subset \Sigma$ of simplices and the set
$B'$ of chambers and will prove in Sections 3 and 4 that the set $B$ is
a basis in $V_\Sigma$ and the set $B'$ is a basis in $V_\Gamma$.
The algorithm of construction of $B$ and $B'$ is inductive on the
dimension $n$ of the affine space $V^n$.

Let an ordering of points $e_1, \ldots , e_N$ satisfy (\ref{conv}). For
each point $e_k$ we construct a set $B_k$ of simplices $\sigma \in
\Sigma$ and a set $B'_k $ of chambers $\gamma \in \Gamma$. Then we
define $B= \bigcup B_k$ and $B' = \bigcup B'_k$.

\smallskip

{\bf First step ($n=2$).} The points $e_1, \ldots , e_N$ lie on the
affine plane $V^2$.
Let us denote by $q$ an edge of a simplex $\sigma \in \Sigma$ and by $Q$
the set of all edges of all simplices $\sigma \in \Sigma$.
Consider the point $e_k \in P_{k-1}$. In $P_{k-1}$ from the point $e_k$
there are following edges $q_i =(e_k, e_i)$, where $i \in (k+1,
\ldots , N)$. Note that since the
points $e \in E$ are not necessarily in general position, some of these
edges may coincide and several points $e_i, \ i\in (k+1, \ldots , N)$
may lie on the same edge.

Let $q_i=(e_k, e_i)$ and $q_j =(e_k, e_j)$, where $i,j \in (k+1, \ldots,
N)$, be two neighbor edges with the vertex $e_k$ (i.e. there is no other
edge $q=(e_k,e_m),\  m \in(k+1, \ldots , N)$ which lies between $q_i$ and
$q_j$).
Let the point $e_i$ be the nearest point of $E$ to the point
$e_k$ on the edge $q_i$ and, respectively, $e_j$ the nearest point to
the point $e_k$ on the edge $q_j$. Let $\sigma= \sigma (e_k, e_i, e_j)$.
We define $B_k$ as the set of all such simplices $\sigma$.
We define then
$$B= \bigcup_k B_k . $$

In Figure 1 there is an example of a set $B$ of simplices constructed
according to this algorithm.

\begin{figure}
\centerline{\psfig{file=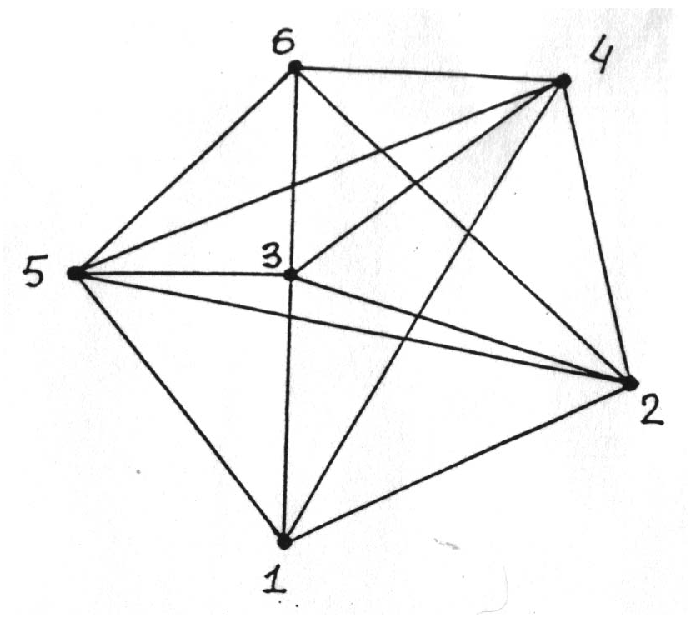,height=5cm}}
\caption{}
\end{figure}

For this example we have: $B_1= (153, 134, 142); \ B_2= (253,236,264);\
B_3=(356,364); \ B_4=(456)$.

\medskip

With the point $e_k$ we also associate the following set of chambers
$B'_k$. In each simplex $\sigma$ let us choose one  chamber\footnote{
In case of an affine plane in each simplex $\sigma$ there is only one
chamber adjacent to the point $e_k$ since every edge of a chamber
$\gamma \in \Gamma$ necessarily lies on some edge $q \in Q$.} adjacent
to the point $e_k$.
We define the set $B'_k$ as the set of all such chambers and define
$$ B' = \bigcup_k B'_k.$$

\begin{figure}
\centerline{\psfig{file=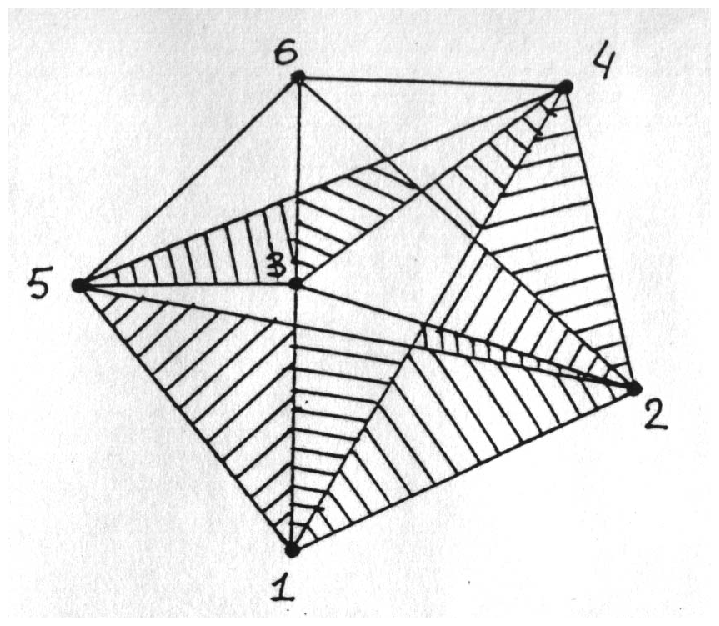,height=5cm}}
\caption{
(There is a mistake in the figure. One
shaded chamber should be in a different place.)}
\end{figure}

In Figure 2 the chambers from the set $B'$ are shaded.

\bigskip

Suppose that we have described the construction
in $V^{n-1}$. Let us describe it in $V^n$.

\smallskip

Let the points $e_i \in E$ be in $V^n$ and let $e_1, \ldots ,
e_N$ be an ordeing satisfying (\ref{conv}). Consider the point $e_k$ and
the hyperplane $H_k$ from Lemma \ref{lemma}. In the hyperplane $H_k$ we
have the set of points $E_k$ (see (\ref{EIK})), the set of simplices
$\Sigma_k$ and the set of chambers $\Gamma_k$. Let us reorder the points
$e^k_i$ so that the ordering $e^k_1, e^k_2, \ldots , e^k _{N_k}$
satisfies (\ref{conv}).

By the induction hypothesis, we can construct a set
of simplices in $H_k$, i.e. the set
$\tilde{B}_k \subset \Sigma_k$ and the set of chambers in $H_k$, i.e.
the set $\tilde{B}'_k \subset \Gamma_k$. Then we define the set $B_k$
of simplices in $V^n$ as $B_k= \{\mu(\sigma^k), \forall
\sigma_k \in \Sigma _k\}$, where the map $\mu$ is defined by
(\ref{mu}). Finally, we define $B=\bigcup_k B_k .$

\smallskip

We have already constructed the set $\tilde{B}'_k$ of chambers in the
hyperplane $H_k$. Consider $\gamma ^k \in \tilde{B}'_k$. According to
the algorithm in $H_k$ the chamber $\gamma^k$ was chosen at a certain
step $j, \ j \in (1, \ldots, N_k)$ and the point $e^k_j$ is a vertex
of $\gamma^k$. Besides, there is one simplex $\sigma^k \in
\tilde{B}_k$ such that $e^k_j$ is a vertex of $\sigma^k$ and $\gamma^k
\subset \sigma^k$. Let us choose a chamber $\gamma \in \Gamma$ such
that

1) $\gamma \subset \mu (\sigma^k) $ ;

2) $\gamma $ is adjacent to
the point $e_k$ and to the edge $(e_k, e_j)$.

Thus, with each chamber $\gamma^k \in \tilde{B}'_k$ we associate a
chamber $\gamma \in \Gamma$.
Denote by $B'_k$ the set of all such chambers $\gamma$ corresponding
to $\gamma^k \in \tilde{B}'_k$ and define $$B'= \bigcup _k B'_k .$$

\smallskip

In Figure 3 there is a fragment of a configuration of points in the
3-dimensional space. In the plane $H_k$, separating the points $e_1,
\ldots, e_k$ and $e_{k+1}, \ldots , e_N$, there are five points which
are reordered according to condition (\ref{conv}). For this ordering
in the plane $H_k$ we construct the basis
$$\tilde{B}_k= \{ (e_{i_1}e_{i_4} e_{i_5}),(e_{i_1}e_{i_5} e_{i_3}),
(e_{i_1}e_{i_3} e_{i_2}); (e_{i_2}e_{i_4} e_{i_5}),
(e_{i_2}e_{i_5} e_{i_3}); (e_{i_3}e_{i_4} e_{i_5})\}$$
of simplices and the basis
$\tilde{B}'_k$ of chambers (are shaded). The set ${B}_k$
consists of the following simplices:
$$B_k= \{ (e_k \bar{e}_{i_1} \bar{e}_{i_4}\bar{e}_{i_5}),\
(e_k \bar{e}_{i_1} \bar{e}_{i_5} \bar{e}_{i_3}), \ (e_k \bar{e}
_{i_1} \bar{e}_{i_3} \bar{e}_{i_2});\
( e_k \bar{e}_{i_2}\bar{e}_{i_4} \bar{e}_{i_5}),\  (e_k \bar{e}_{i_2}
\bar{e}_{i_5} \bar{e}_{i_3});\
( e_k \bar{e}_{i_3} \bar{e}_{i_4}\bar{e}_{i_5})\},$$
where $\bar{e}_j= \mu (e^k_j).$

\begin{figure}
\centerline{\psfig{file=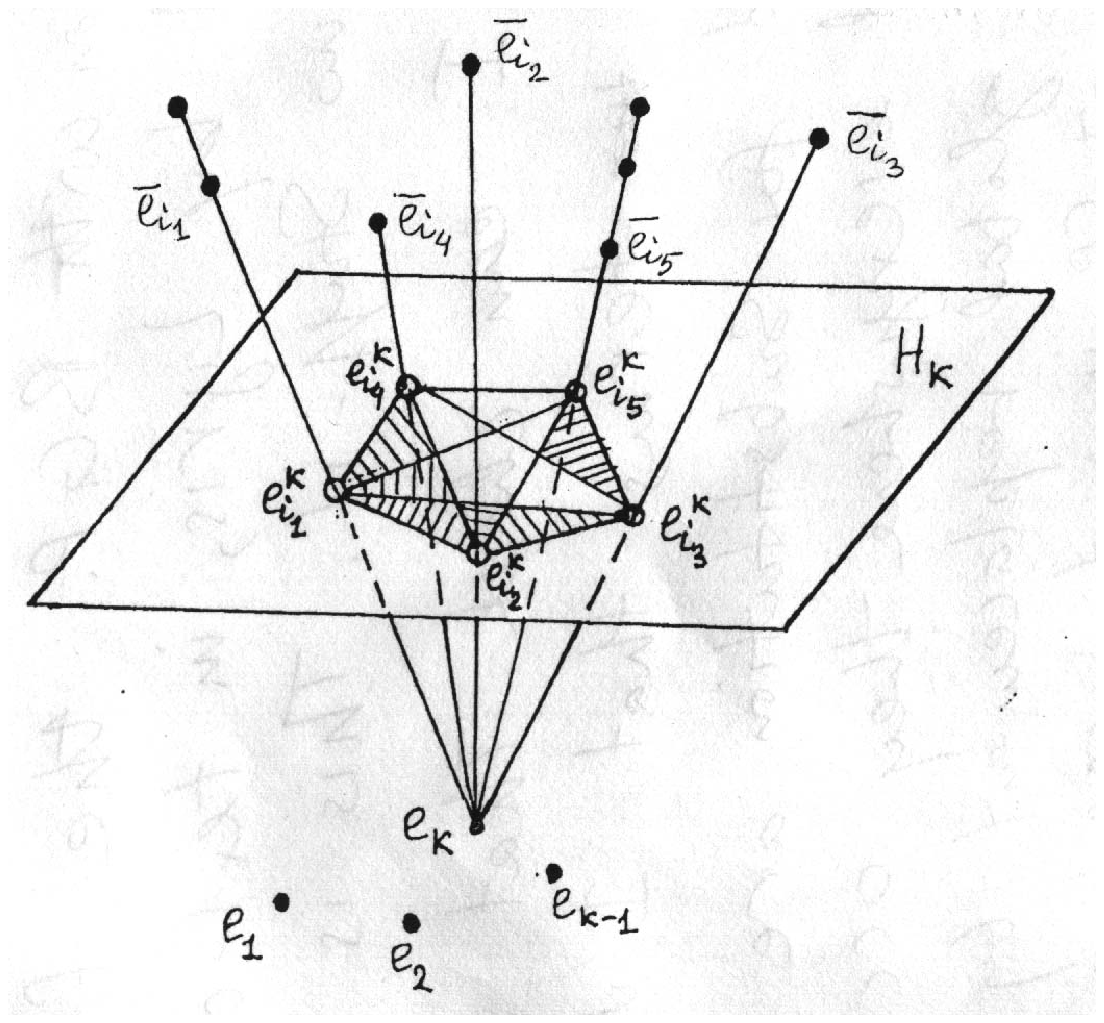,height=9cm}}
\caption{}
\end{figure}

Note that a chamber from $B'_k$ cannot be seen in Figure 3 since we need
to take into account the points $e_1, \ldots, e_{k-1}$ which are below
the plane $H_k$.

\section{ Linear independence of simplices $\sigma \in B$ and linear
independence of chambers $\gamma \in B'$.}

Let $B \subset \Sigma$ and $B' \subset \Gamma$ be the sets constructed
in Section 2.

\begin{theorem}
\label{indep}
The simplices $\sigma \in B$ are linearly independent in $V_\Sigma$ and
the chambers $\gamma \in B'$ are linearly independent in $V_\Gamma$.
\end{theorem}

First let us prove the following proposition.

\begin{proposition}
\label{block}
Let $\tilde{A}$ be a submatrix of the incidence matrix $A$
corresponding to the rows $\sigma \in B$ and the columns $\gamma \in
B'$. The columns and rows of the matrix $\tilde{A}$ can be ordered in
such a way that: 1) $\tilde{A}$ is a block matrix $\parallel
\cal{A}_{ik} \parallel $, where $\cal{A}_{ik}= 0$ for $i>k$, and 2) each
 diagonal element $a_{\sigma, \gamma}$ of the matrix $\tilde{A}$ equals
 $1$, i.e. $a_{\sigma, \gamma}=1.$

\end{proposition}

{\bf Proof.} Let us recall that by the algorithm the set $B$ of
simplices was constructed as $B=\cup B_k$ and the set $B'$ of chambers
as $B'=\cup B'_k$, where $B_k, B'_k$ correspond to the point $e_k$ in
the ordering $e_1, \ldots, e_N$.
Thus, the matrix $\tilde{A}$ is a block-matrix $\parallel \cal{A}_{ik}
\parallel$, where

\begin{equation}
\label{aik}
\cal{A}_{ik}= \{ a_{\sigma, \gamma}, \ \sigma \in
B_i, \ \gamma \in B'_k\}.
\end{equation}

Consider a diagonal block $\cal{A}_{kk}$ ($k=1, \ldots, N-n$).
By the construction, to each simplex $\sigma \in B_k$ there corresponds
a chamber $\gamma \in B'_k$ such that $\gamma \subseteq \sigma$,
therefore, $a_{\sigma, \gamma}=1$. Thus, if we choose the corresponding
orderings of columns in $B'_k$ and rows in $B_k$, we obtain $"1"$ on the
main diagonal in the block $\cal{A}_{kk}$ and, therefore, any diagonal
element $a_{\sigma, \gamma} $ of the matrix $\tilde{A}$ is such that
$a_{\sigma, \gamma} =1$.

\smallskip

Consider a block $\cal{A}_{ik}$, where $i>k$.
We need to prove that
$a_{\sigma, \gamma}= 0$,  where $\sigma \in B_i$ and $\gamma \in B'_k$.
Let $\sigma \in B_i$ and $\gamma \in B'_1 \cup \ldots \cup B'_{i-1}$.
Then $\sigma \in P_{i-1}=conv(e_i, \ldots, e_N)$.
>From the algorithm for the construction of the set $B'$ it is easy to
see that $\gamma \in S_1 \cup \ldots \cup S_{i-1}$. Due to formulas
(\ref{decomp}) and (\ref{SP}) we conclude that $\gamma \not \subset
\sigma$ and, therefore, $a_{\sigma, \gamma}=0$. Thus, the matrix
$\tilde{A}$ is upper triangular as a block matrix. \qed

\medskip

{\bf Proof of Theorem \ref{indep}.}  We will show that the submatrix
$\tilde{A}$ defined above is an upper triangular matrix.
Due to Proposition \ref{block} it is sufficient to prove that a block
$\cal{A}_{kk}$ ( where $k=1, \ldots, N-n$) of the matrix $\tilde{A}$ is
upper triangular.

\smallskip

Let us consider the case {\bf $n=2$}. Let $e_1, \ldots, e_N$ be an
ordering satisfying (\ref{conv}). By the construction, the
simplices $\sigma \in B_k$ lie between the neighbor edges starting at
the point $e_k$, therefore, the open simplices $\sigma \in B_k$ are
disjoint. Besides, there is exactly one chamber $\gamma \in B'_k$ (i.e.
with the vertex $e_k$) such that $\gamma \in \sigma$, therefore, the
block $\cal{A}_{kk}$ is the identity matrix.

\medskip
Since for a general $n$ the notations are cumbersome
we consider in detail the case $n=3$ which already contains all the
technique.
Let $e_1, \ldots, e_N$ be an ordering
satisfying (\ref{conv}). Consider a point $e_k$ and the block
$\cal{A}_{kk}$ of the matrix $\tilde{A}$. We recall (see Section 2.3)
that a simplex $\sigma \in B_k$ is defined as $\mu (\sigma^k)$, where
the 2-dimensional simplex $\sigma ^k \in \tilde{B}_k$ lies in the plane
$H_k$ separating the points $e_1, \ldots, e_k$ and $e_{k+1}, \ldots,
e_N$.

The set $\tilde{B}_k$ of simplices is constructed according to the
algorithm in $H_k$. For this we order the points $e^k_i $ (see
(\ref{EIK})) according to condition (\ref{conv}).
Let $e^k_1, \ldots, e^k_{N_k}$ be such an ordering.
Then $\tilde{B}_k= \bigcup _i \tilde{B}_{k,i}$, where $\tilde{B}_{k,i}$
is the set of 2-dimensional simplices which were chosen in the algorithm
for the point $e^k_i$. We have a similar equality for chambers:
$\tilde{B}'_k=\bigcup_i \tilde{B}'_{k,i}$, where $\tilde{B}'_{k,i}$ is
the set of 2-dimensional chambers chosen in the algorithm for the point
$e^k_i$.

Due to the inductive construction of simplices $\sigma \in B_k$ and
chambers $\gamma \in B'_k$ we obtain also the following formulas:

$$B_k= \bigcup _i B_{k,i}$$
 and $$B'_k=\bigcup_i B'_{k,i}.$$
This means that the block $\cal{A}_{kk}$ consists, in turn, of blocks
$\cal{B}_{i_1, k_1}$:

$$\cal{B}_{i_1,k_1}=\{a_{\sigma, \gamma}: \sigma^k \in \tilde{B}_{k,i_1},
\gamma^k \in \tilde{B}'_{k,k_1}\}, $$ where
$\sigma = \mu(\sigma^k)$
and the chamber $\gamma$ corresponds \footnote{Note that in the
algorithm there is no direct correspondence between the chambers
$\gamma$ and $\gamma^k$; given a chamber $\gamma^k \subset H_k$ we find
the corresponding simplex $\sigma^k \subset H_k$, then in the simplex
$\sigma =\mu(\sigma^k)$ we choose a certain chamber (see Section 2.3).}
to the chamber $\gamma^k \subset H_k$.

\smallskip

$1^0$. Consider a diagonal block $\cal{B}_{k_1,k_1}$ of the block
$\cal {A}_{kk}$, i.e. all $a_{\sigma, \gamma}$, where $\sigma \in
B_{k,k_1}$ and $\gamma \in B'_{k,k_1}$.

Let $\sigma, \sigma_0 \in B_{k,k_1}$. Then the corresponding simplices
$\sigma^k $ and $\sigma ^k_0$ were constructed on the plane $H_k$ from
the same point $e^k_{k_1}$. By the algorithm the open simplices
$\sigma^k,\sigma ^k_0$ are disjoint. Then (see Proposition
\ref{prop1}) $\mu(\sigma^k) \cap \mu(\sigma^k_0) =\emptyset$, where $\mu
$ is defined by formula (\ref{mu}). But
$\sigma=\mu(\sigma^k)$, and $\sigma_0 =\mu(\sigma^k_0)$. Obviously, for
any chamber $\gamma$  such that $\gamma \subset \sigma$, we have $\gamma
\not \subset \sigma_0$. This means that the diagonal block
$\cal{B}_{k_1,k_1}$ is the identity matrix.

\medskip

$2^0$. Let us show that $\cal{B}_{i_1,k_1}= 0$ for $i_1>k_1$.
Indeed, let $\sigma \in B_{k,i_1}$ and $\gamma \in
B'_{k,k_1}$. Then from the construction
we have $\sigma^k \in \tilde{B}_{k,i_1}$ and
$\sigma^k \in P^k_{i_1-1}=conv(e^k_{i_1}, \ldots, e^k_{N_k})$ (see
(\ref{PS})).

Concerning $\gamma$ we know that there is a simplex $\sigma_0 \in
B_{k,k_1}$ such that:

1) $\gamma \in \sigma_0$;

2) $\gamma$ is adjacent to the point $e_k$;

3) $\gamma$ is adjacent to the edge $(e_k, e^k_{k_1})$.

Since $\sigma_0 \in B_{k,k_1}$ we have $\sigma^k_0 \in
\tilde{B}_{k,k_1}$ and $\sigma^k_0 \subset  P^k_{k_1-1}$.

For $i_1>k_1$ we can rewrite the formula (\ref{PK})
as follows:

$$P^k= S^k_1 \cup \ldots \cup S^k_{k_1} \cup S^k_{k_1+1} \cup
\ldots \cup S^k_{i_1-1}\cup P^k_{i_1-1},$$
where all the open polytopes $S^k_j$ and $P^k_{i_1-1}$ are disjoint.
Note that if $i_1=k_1+1$ then $S^k_{k_1}=S^k_{i_1-1}$.

Due to decomposition (\ref{tau}) we obtain
$$\tau(P^k)= \tau(S^k_1) \cup \ldots \cup \tau(S^k_{k_1}) \cup
\tau(S^k_{k_1+1}) \cup
\ldots \cup \tau(S^k_{i_1-1})\cup \tau(P^k_{i_1-1}),$$ where all the
corresponding open cones are disjoint.
We have $\sigma= \mu(\sigma^k) \subset \tau(P^k_{i_1-1})$.
\smallskip
Let us show that $\gamma \subset \tau(S^k_{k_1})$.
Since $\gamma \subset \sigma_0$, then $\gamma \subset \mu(\sigma^k_0)$,
i.e. $\gamma \subset \ \tau(P^k_{k_1-1})$.
Note that
$P^k_{k_1-1}= S^k_{k_1} \cup S^k_{k_1+1} \cup \ldots$ (see formulas
(\ref{PS}) and (\ref{PK})).
Then $\gamma \subset \tau(S^k_{k_1}) \cup \tau(S^k_{k_1+1}) \cup
\ldots.$ Any chamber $\gamma \in \Gamma$ may lie in only one of these
cones.
Since $\gamma $ is adjacent to the edge $(e_k, e^k_{k_1})$ we
obtain  $\gamma \subset \tau(S^k_{k_1})$.
Then $\gamma \not \subset \tau(P^k_{i_1-1})$, i.e. $\gamma \not \subset
\sigma$. We have proved that $a_{\sigma, \gamma}=0$ for $\sigma \in
B_{k,i_1}$ and $\gamma \in B'_{k,k_1}$ if $i_1>k_1$.

This completes the  proof (for $n=3$) that the matrix $\tilde{A}$ is an
upper triangular matrix.

\bigskip

In case of $n$-dimensional space an element $a_{\sigma, \gamma}$ of the
matrix $\tilde{A}$ belongs to a sequence of enclosed blocks which can
be denoted as $\cal{A}_{i_1,k_1}, \cal{B}_{i_2,k_2}, \cal{C}_{i_3,k_3},
\ldots$. (One can check that the "depth" of the enclosed blocks, i.e.
the length of the sequence, is $n-1$.)

$1^0.$ Let us show that a diagonal block of maximal depth is the
identity matrix. Let $\cal{Z}_{k_{n-2}, k_{n-2}}$ be such a block. Then
$\cal{Z}_{k_{n-2}, k_{n-2}}$ lies inside all diagonal blocks
$\cal{A}_{k_0,k_0}, \cal{B}_{k_1,k_1}, \cal{C}_{k_2,k_2}, \ldots$ of the
matrix $\tilde{A}$. Let $\sigma, \sigma_0$ be two simplices from the
block $\cal{Z}_{k_{n-2},k_{n-2}}$. This  means that $\sigma$ and
$\sigma_0$ have the same sequence of corresponding points in the
inductive construction.
Let us denote these points by $e_{k_0},
e^{(1)}_{k_1}, e^{(2)}_{k_2}, \ldots, e^{(n-2)}_{k_{n-2}}$, where
$e_{k_0}\in E$, $e^{(1)}_{k_1} \in H_{k_0}$, $e^{(2)}_{k_2}$ lies in the
$(n-2)$-dimensional plane which corresponded to the point
$e^{(1)}_{k_1}$ in the algorithm, and so on; finally,
the point $e^{(n-2)}_{k_{n-2}}$ lies in the
$2$-dimensional plane which corresponded to the point
$e^{(n-3)}_{k_{n-3}}$. On this latter plane there are two simplices
$\tilde{\sigma}, \tilde{\sigma_0}$ which were both chosen for the point
$e^{(n-2)}_{k_{n-2}}$.

According to the algorithm for the plane the open simplices
$\tilde{\sigma}, \tilde{\sigma_0}$ are disjoint. Therefore (see
Proposition \ref{prop1}), for the $3$-dimensional open simplices we have
$\mu(\sigma) \cap \mu(\sigma_0) =\emptyset$, where $\mu$ is a map
defined by formula (\ref{mu}) for the point
$e^{(n-3)}_{k_{n-3}}$ instead of $e_k$.
Then in order to obtain $4$-dimensional simplices, the map $\mu$ for the
point $e^{(n-4)}_{k_{n-4}}$ was applied to the simplices
$\mu(\sigma), \mu(\sigma_0)$ and so on. Finally, the map $\mu$ for the
point $ e_{k_0}$ gives the simplices $\sigma , \sigma_0 \in \Sigma$.
Clearly, we have $\sigma \cap \sigma_0 =\emptyset$. Therefore, for any
chamber $\gamma \subset \sigma$ we have $\gamma \not \subset \sigma_0$,
i.e. the block $\cal{Z}_{k_{n-2}, k_{n-2}}$ is the identity matrix.

\medskip

$2^0.$ In order to prove that the matrix $\tilde{A}$ is upper
triangular it suffices to check that  for $i>j$ any block
$\cal{X}_{i,j} \in (\cal{B}_{i_1,k_1}, \cal{C}_{i_2,k_2}, \ldots) $
satisfies the condition $\cal{X}_{i,j}=0$.
In Proposition \ref{block} we have proved that for $i>k$ a block
$\cal{A}_{ik}=0$.
By applying the arguments from $2^0$ of $n=3$ case, we can prove that in
a diagonal block $\cal{A}_{kk}$ we have $\cal{B}_{{i_1},{k_1}}=0$ for
$i_1>k_1$. Similarly, $\cal{C}_{{i_2},{k_2}}=0$ for $i_2>k_2$, etc.

It follows from $1^0$ and $2^0$ that  the submatrix $\tilde{A}$ of the
incidence matrix $A$ is an upper triangular matrix. \qed

\section{ The set $B$ is a geometrical basis of simplices.}

In this section we will prove that the set $B$ constructed by the
algorithm is a geometrical basis. First we
repeat some of the definitions and theorems from \cite{AGZ} and
\cite{A}.

Let $E=(e_1, \ldots, e_N)$ be a set of points in the $n$-dimensional
affine space $V^n$, $\Sigma$ the set of simplices  and $\Gamma$ the set
of chambers (defined in the Introduction).

Consider a subset $S \subseteq E$ consisting of $n+2$ points and such
that $S$ contains at least $n+1$ points in general position.
Denote
\begin{equation}
\label{fsimpl}
f = \{ \sigma: \sigma(e_{i_1}, \ldots, e_{i_{n+1}}) \in \Sigma \hbox{
and } e_{i_k} \in S \}
\end{equation}
Thus, with each $S$ we associate a subset $f\subset \Sigma$
(clearly, $f \neq \emptyset$). Let $F$ be the set of all such $f$
corresponding to  all possible $S \subseteq E$.

\begin{definition}
We say that a point $p$ of affine space is visible from a point
$e$, $e \neq p$, with respect to a simplex $\sigma$ if the open
segment $(e,p)$ ia disjoint with $\sigma$, i.e. $(e,p) \cap \sigma
=\emptyset $.

We say also that a subset $S$ of points is visible from the point $e$ if
every point of this subset is visible from $e$.
\end{definition}

\begin{theorem}
\footnote{This theorem is stated in \cite{AGZ} in
another form. Here we use the important geometric notion of
visibility.}
\label{TS}
Let $\sigma \in \Sigma$ be a simplex and $e\in E $ a
point that is not a vertex of $\sigma$.

There is the following linear relation in $V_\Sigma$ among
simplices:

\begin{equation}
\label{RelS}
\sigma = \sum _{q_i \in Q^+} \sigma (q_i,e) - \sum _{q_i \in Q^-}
\sigma (q_i,e),
\end{equation}

where

$Q^+$ is the set of all facets (i.e. ($n-1$)-dimensional faces )  $q_i$
of the simplex  $\sigma$ that are not visible from $e$ (with respect to
$\sigma$);

$Q^-$ is the set of all facets $q_i$ of the simplex  $\sigma$ that are
visible from $e$ (with respect to $\sigma$);

$\sigma (q_i,e)$ is the $n$-dimensional simplex spanned by the facet
$q_i$ of the simplex $\sigma$ and the point $e$.
\end{theorem}

\begin{definition}
Let $B \subset \Sigma$. We say that an element $\sigma  \not \in B$ is
expressed in one step ``in terms of the set $B$ using $F$'' if there
exists $f\in F$ such that $\sigma \in f, \ and \ f \setminus \sigma
 \subseteq B$.

We say that an element $\sigma  \not \in B$ can be expressed in $k$
steps in terms of the set $B$ using $F$ if there exists a sequence
$\sigma_1, \ldots , \sigma_k,\  \sigma_i \in \Sigma,$ such that
$\sigma_k=\sigma$ and $\sigma_1$ is expressed in
one step in terms of the set $B$ (using $F$), $\sigma_2$ is expressed in
one step in terms of the set $B \cup \sigma_1, \ldots ,
\ \sigma_k$ is expressed in one step in terms of the set $B\cup
\sigma_1 \ldots \cup \sigma_{k-1}$.

\end{definition}

\begin{definition}
A subset $B \subset \Sigma $ is a geometrical basis in $V_\Sigma$
 with respect to $F$ if it satisfies the following conditions:

1) $B$ is a basis in $V_\Sigma$;

2) for any $\sigma \in \Sigma, \ \sigma  \not \in B$ there exists $k$
such that $\sigma $ can be expressed in $k$ steps in terms of $B$ using
$F$.
\end{definition}

Let $B \subset \Sigma$ be the set of simplices constructed by
the algorithm of Section 1.

\begin{theorem}
\label{TS1}
The set $B$ is a geometrical basis in $V_\Sigma$ with respect to the
system $F$ defined by formula (\ref{fsimpl}).
\end{theorem}

{\bf Proof.} Due to Theorem \ref{indep} the simplices $\sigma \in B$ are
linearly independent, therefore it is sufficient to prove that any
simplex $\sigma \in \Sigma$ can be expressed in a finite number of steps
in terms of $B$ with respect to the system $F$. This will be proved by
induction on the dimension $n$ of the space $V^n$.

{\bf Induction on $n$.} {\em First step ($n=2$).} The theorem is proved
in \cite{AG}.
\smallskip

{\em Passing from $n-1$ to $n$.}
Suppose that the statement is true for $V^{n-1}$. Let us prove it for
$V^n$.

Let the points $e_i \in E$ be in $V^n$ and an ordering $e_1, \ldots ,
e_N$ satisfy (\ref{conv}). Let us denote by $\Sigma^i$
the set of simplices $\sigma \in \Sigma$ such that $\sigma$ has the
point $e_i$ as the vertex with the minimal number. We have
\begin{equation} \label{sigma} \Sigma= \Sigma^1 \cup \ldots \cup
\Sigma^{N-n}= \bigcup_{i=n}^{N-1} \Sigma^{N-i},
\end{equation} where
$\Sigma^i \cap \Sigma^j=\emptyset,\ i\neq j$. We also have $\Sigma^N=
\Sigma^{N-1}= \ldots = \Sigma^{N-n+1}=\emptyset$. The set $\Sigma^{N-n}$
either contains only one simplex $\sigma (e_N,e_{N-1}, \ldots, e_{N-n})$
or $\Sigma^{N-n}=\emptyset$.

We need to prove that any simplex $\sigma \in \Sigma$ can be expressed
in terms of $B$ using $F$. Due to the partition
(\ref{sigma}) we can prove this statement by induction on $i$
considering cases when $\sigma \in \Sigma^{N-i}$, where $i=n, \ldots
,N-1$.

{\bf Induction on $i$.} {\em First step.}
Let us check the first nontrivial step of induction. Let
$\Sigma^{N-i_0}$ (where $i_0  \geq n$) be the first nonempty set, i.e.
$$\Sigma^{N-i_0}=\{ \sigma(e_N, \ldots , e_{N-i_0+1}, e_{N-i_0}) \},$$
where the points $e_N, \ldots , e_{N-i_0+1}$ lie in an
$(n-1)$-dimensional hyperplane, while the point $e_{N-i_0}$ does not lie
in this hyperplane. It is clear from the algorithm that the simplex
$\sigma(e_N, \ldots , e_{N-i_0+1}, e_{N-i_0})$ belongs to the set
$B_{N-{i_0}} \subset B$
and, therefore, can be expressed in terms of $B$ in $0$ steps.

\medskip

{\em Passing from $i-1$ to $i$.}
Suppose that any simplex $\sigma \in \Sigma^{N-(i-1)}$ can be
expressed in terms of $B$ using $F$.
Let us show that any simplex $\sigma_0 \in \Sigma^{N-i}$ can also be
expressed in terms of $B$ using $F$.
\medskip

Consider the point $e_{N-i}$ and the hyperplane $H=H_{N-i}$ of the
algorithm.
According to the algorithm we mark in $H$ the set $E_{N-i}$ of points
$e_1^{N-i}, e_2^{N-i}, \ldots, e_{N_i}^{N-i}$, where
$e_j^{N-i}=\overrightarrow{(e_{N-i},e_j)} \cap H$.
To any vertex of $\sigma_0$ other than $e_{N-i}$, there corresponds a
point from $E_{N-i}$, therefore, to the simplex $\sigma_0$ there
corresponds one $(n-1)$-dimensional simplex $\sigma_0^{N-i}$ in $H$.

Let $\Sigma_{N-i}$ be the set of all ($n-1$)-dimensional simplices with
the vertices in $E_{N-i}$. We have $\sigma_0^{N-i} \in \Sigma_{N-i}$.
Applying the algorithm in the hyperplane $H$, we construct the set
$\tilde{B}_{N-i}$ of simplices, $\tilde{B}_{N-i} \subset \Sigma_{N-i}$.
Due to the assumption of the induction on $n$, the simplex
$\sigma_0^{N-i} \in \Sigma_{N-i}$ can be expressed in terms of
$\tilde{B}_{N-i}$ using $F'$ (where $F'$ is the set of all subsets $f'$
defined in $H$ by formula (\ref{fsimpl}) ).

\smallskip

First, let us show that the simplex $\mu(\sigma_0^{N-i})$, where
$\mu$ is the map defined by (\ref{mu}) for $k=N-i$,
can be expressed in terms of $B$ using $F$.

As we have proved, the simplex $\sigma^{N-i}_0$ can be expressed in
terms of $\tilde{B}_{N-i}$ using $F'$ in a finite number of steps, for
example, in $k$ steps.
This means that there exists a sequence of simplices
$\sigma_1^{N-i}, \ldots , \sigma_k^{N-i}$, where $\sigma_k^{N-i}=
\sigma^{N-i}_0$ and a sequence of $f_1',  \ldots , f_k ',\ \  f_i' \in
F'$, such that

\begin{equation}
\label{first}
f_1' \setminus \sigma_1^{N-i} \subseteq \tilde{B}_{N-i},
\end{equation}
\begin{equation}
\label{second}
f_2' \setminus \sigma_2^{N-i} \subseteq
(\tilde{B}_{N-i} \cup \sigma_1^{N-i}),
\end{equation}
$$ \ldots  $$
\begin{equation}
\label{last}
f_k' \setminus \sigma_k^{N-i} \subseteq (\tilde{B}_{N-i} \cup
\sigma_1^{N-i} \cup \ldots \cup \sigma_{k-1}^{N-i}).
\end{equation}

Consider formula (\ref{first}). It means that the simplex
$\sigma_1^{N-i}$ can be expressed in terms of $\tilde{B}_{N-i}$ using
$F'$. Let us show that this implies that the simplex
$\mu(\sigma_1^{N-i})$ can be expressed in terms of $B_{N-i}$ using $F$.
Indeed, according to the definition, the set $f'_1$
contains all the simplices with vertices in some $n+1$ points
$e^{N-i} _{i_1}, \ldots, e^{N-i}_{i_{n+1}} \in E_{N-i}$.
One of these simplices is
$\sigma_1^{N-i}$. For simplicity of notations let us assume that
$\sigma_1^{N-i}=\sigma( \widehat{e^{N-i}_{i_1}}, e^{N-i}_{i_2}, \ldots,
e^{N-i}_{i_{n+1}}), $
where  $\widehat{e}$ means that the point $e$ is not a vertex of the
simplex.
Thus, $f_1'$ consists of the following simplices:

$$\sigma_1^{N-i}=\sigma( \widehat{e^{N-i}_{i_1}}, e^{N-i}_{i_2}, \ldots,
e^{N-i}_{i_{n+1}}) ,$$
$$\sigma( e^{N-i}_{i_1}, \widehat{e^{N-i}_{i_2}}, \ldots,
e^{N-i}_{i_{n+1}}) ,$$
$$\ldots$$
$$\sigma( e^{N-i}_{i_1}, e^{N-i}_{i_2}, \ldots,
\widehat{e^{N-i}_{i_{n+1}}}), $$
where each simplex except $\sigma_1^{N-i}$ belongs to $\tilde{B}_{N-i}$
due to formula (\ref{first}).

\smallskip

To $f'_1 \in F'$ let us associate an element $f_1 \in F$. For this
consider $n+2$ points $e_{N-i}, e^\mu _{i_1}, \ldots,
e^\mu_{i_{n+1}} \in E$, where $e_j^\mu$ is the
point of $E$ which is the closest to the point
the point $e_{N-i}$ among all points of $E$ lying on the ray
$\overrightarrow{(e_{N-i},e^{N-i}_j)}$. Then $f_1$ is the set of all
simplices with the vertices in these $n+2$ points. Thus, $f_1$ consists
of the following simplices:

$$\sigma(e^\mu_{i_1}, e^\mu_{i_2}, \ldots, e^\mu_{i_{n+1}}) ,$$
$$\sigma(e_{N-i}, \widehat{e^\mu_{i_1}}, e^\mu_{i_2}, \ldots,
e^\mu_{i_{n+1}}) ,$$
$$\sigma(e_{N-i}, e^\mu_{i_1}, \widehat{e^\mu_{i_2}}, \ldots,
e^\mu_{i_{n+1}}) ,$$
$$\ldots$$
$$\sigma(e_{N-i}, e^\mu_{i_1}, e^\mu_{i_2}, \ldots,
\widehat{e^\mu_{i_{n+1}}}). $$

The point $e_{N-i}$ is not the vertex of the
simplex $\sigma(e^\mu _{i_1}, \ldots, e^\mu_{i_{n+1}}) $. Therefore,
$\sigma(e^\mu _{i_1}, \ldots, e^\mu_{i_{n+1}}) \in \Sigma^{N-(i-1)}$.
By induction hypothesis for $i$ this simplex can be
expressed in terms of $B$ using $F$.

Note that according to the algorithm,
$$\sigma (e_{N-i}, e^\mu _{i_1}, \ldots, \widehat{e^\mu_{i_j}},
\ldots, e^\mu_{i_{n+1}}) =\mu (\sigma(e^{N-i}_{i_1}, \ldots,
\widehat{e^{N-i}_{i_j}}, \ldots, e^{N-i}_{i_{n+1}}))$$ for any $j =1,
\ldots, n+1$. Therefore, we have
$$\sigma(e_{N-i}, \widehat{e^\mu_{i_1}}, e^\mu_{i_2}, \ldots,
e^\mu_{i_{n+1}}) =\mu (\sigma^{N-i}_1).$$

Since each simplex $\sigma(e^{N-i}_{i_1}, \ldots,
\widehat{e^{N-i}_{i_j}}, \ldots, e^{N-i}_{i_{n+1}}), \ j =2,
\ldots, n+1$, belongs to $\tilde{B}_{N-i}$, each simplex
$\sigma (e_{N-i}, e^\mu _{i_1}, \ldots, \widehat{e^\mu_{i_j}},
\ldots, e^\mu_{i_{n+1}})$ belongs to $B_{N-i}$ according to the
construction of $B_{N-i}$.
Thus,
$f_1 \setminus \mu(\sigma_1^{N-i}) \subseteq (B_{N-i} \cup
\sigma(e^\mu _{i_1}, \ldots, e^\mu_{i_{n+1}}),$ where, as we have
already mentioned, the simplex $\sigma(e^\mu _{i_1}, \ldots,
e^\mu_{i_{n+1}})$ can be expressed in terms of $B$ using $F$.
We have proved that $\mu(\sigma_1^{N-i})$ can be expressed in terms of
$B$ using $F$.
\smallskip

Considering consecutively formulas (\ref{second}) -- (\ref{last}) we can
prove similarly that the simplex $\mu(\sigma_0^{N-i})$ can be expressed
in terms of $B$ using $F$.

\medskip

Consider the simplex $\sigma_0$.
The following cases are possible:

1) $\sigma_0= \mu(\sigma_0^{N-i})$; in this case the proof is already
finished.

2) $\sigma_0 \neq \mu(\sigma_0^{N-i})$.
Let $\sigma_0^{N-i}= \sigma(e_{j_1}^{N-i}, \ldots, e_{j_n}^{N-i})$.
This case can only occur if there is a vertex of $\sigma_0$ which lies
on some ray $\overrightarrow{(e_{N-i}, e_{j_k}^{N-i})},$ where $k\in
(1, \ldots, n)$, and which is not the closest point from $E$ to the
point $e_{N-i}$.

Let us show that since the simplex $\mu (\sigma_0^{N-i})$ can be
expressed in terms of $B$ using $F$ then the simplex $\sigma_0$ can be
also expressed in terms of $B$ using $F$.

For simplicity of notations let us denote here vertices of the
simplices $\sigma_0$ and $\mu(\sigma_0^{N-i})$ as follows:
$$\sigma_0=\sigma(e_{N-i},e_1^0, \ldots, e_n^0),$$
$$\mu(\sigma_0^{N-i})=\sigma(e_{N-i}, e_1^\mu, \ldots , e_n ^\mu).$$
We can also assume that in this notation the points $e_k^0$ and
$e_k^\mu$ lie on the same ray $\overrightarrow{(e_{N-i},
e_{j_k}^{N-i})}$.

First, let us consider the case when the simplices $\sigma_0$ and
$\mu(\sigma_0^{N-i})$ differ only by one point, i.e.
$e_k^0 \neq e_k^ \mu$ for some $k$ and we have
$$\sigma_0= \sigma(e_{N-i}, e_1^\mu, \ldots, e_{k-1}^\mu, e_k^0,
e_{k+1}^\mu, \ldots, e_n^\mu),$$
$$\mu(\sigma_0^{N-i})=\sigma(e_{N-i}, e_1^\mu, \ldots, e_{k-1}^\mu,
e_k^\mu, e_{k+1}^\mu, \ldots, e_n^\mu).$$

Then for the $n+2$ points $e_{N-i}, e_1^\mu, \ldots,
e_{k-1}^\mu, e_k^\mu, e_{k+1}^\mu, \ldots, e_n^\mu, e_k^0 \in E$
there exists $f \in F$ which consists of the following simplices:

$$\sigma( e_{N-i}, e_1^\mu, \ldots, e_{k-1}^\mu, e_k^\mu, e_{k+1}^\mu,
\ldots, e_n^\mu)=\mu(\sigma_0^{N-i}),$$
$$\sigma(e_{N-i}, e_1^\mu, \ldots, e_{k-1}^\mu, e_k^\mu, e_{k+1}^\mu,
\ldots,e^\mu_{n-1}, e_k^0),$$
$$\ldots$$
$$\sigma(e_{N-i}, e_1^\mu, \ldots, e_{k-1}^\mu, e_k^\mu, e_k^0,
\ldots, e_n^\mu)$$
$$\sigma( e_{N-i}, e_1^\mu, \ldots, e_{k-1}^\mu, e_k^0, e_{k+1}^\mu,
\ldots, e_n^\mu)=\sigma_0,$$
$$\sigma(e_{N-i}, e_1^\mu, \ldots, e_k^0, e_k^\mu, e_{k+1}^\mu,
\ldots, e_n^\mu),$$
$$\ldots$$
$$\sigma(e_{N-i}, e_k^0, \ldots, e_{k-1}^\mu, e_k^\mu, e_{k+1}^\mu,
\ldots, e_n^\mu),$$
$$\sigma(e_k^0, e_1^\mu, \ldots, e_{k-1}^\mu, e_k^\mu, e_{k+1}^\mu,
\ldots, e_n^\mu).$$

Since the points $e_{N-i}, e_k^0, e_k^\mu$ lie on the same ray
$\overrightarrow{(e_{N-i}, e_{j_k})}$, all the simplices except
$\sigma_0$ and $\mu(\sigma_0^{N-i})$ are not $n$-dimensional and do not
belong to $\Sigma$. This implies that the simplex $\sigma_0$ can be
expressed in terms of $B$ using $f \in F$.

If the simplices $\sigma_0$ and $\mu(\sigma_0^{N-i})$ differ by two
points, then we consider first a simplex $\sigma'$ which differs from
the simplex $\mu(\sigma_0^{N-i})$ only by one point and express it in
terms of $B$ using $F$. Then by applying the above arguments to
the simplices $\sigma'$ and $\sigma_0$ we can express the simplex
$\sigma _0$ in terms of $B$ using $F$. By repeating this process we
prove that the simplex $\sigma_0$ can be expressed in terms of $B$ using
$F$. Thus, the set $B$ constructed by the algorithm is a
geometrical basis in $V_\Sigma$. \qed

\smallskip

A pair of a basis $e_1, \ldots, e_n$ in the space $V$ and a basis
$f_1, \ldots, f_n$ in the dual space $V'$ is called a {\em triangular
pair} if $(e_i, f_k)=0$ for $i>k$ and $(e_i,f_i)=1$.

\smallskip

Theorems \ref{TS1} and \ref{indep} imply

\begin{theorem}
The set $B'\subset\Gamma$ constructed by the algorithm is a basis in
$V_\Gamma$. The basis of simplices $B$ and the basis of chambers $B'$
form a triangular pair.
\end{theorem}

\end{document}